\documentclass[12pt]{article}
\usepackage{amsfonts}
\usepackage{latexsym}
\begin{document}
 
\sloppy
 
\newcommand{\e}{\epsilon}
\renewcommand{\a}{\alpha}
\renewcommand{\b}{\beta}
\newcommand{\om}{\omega}
\newcommand{\Om}{\Omega}
\newcommand{\p}{\partial}
\renewcommand{\phi}{\varphi}
\newcommand{\N}{{\mathbb N}}
\newcommand{\R}{{\mathbb R}}
\newcommand{\cF}{{\cal F}}

\newtheorem{theorem}{Theorem}
\newtheorem{lemma}{Lemma}

\title{Dissipative Quasigeostrophic Dynamics under Random Forcing
\footnote{This work was supported by the National Science
                Foundation Grant DMS-9704345.}    } 

\author{  }
 
\date{April 6, 1998    }
\maketitle
 
\begin{center}
 
James R. Brannan$^{1}$, Jinqiao Duan$^{1}$, and Thomas Wanner$^{2}$  \\

\bigskip

{\em 
1. Department of Mathematical Sciences, Clemson University,\\
   Clemson, South Carolina 29634, USA. \\
   E-mail: duan@math.clemson.edu \\
   
2. Institute of Mathematics, University of Augsburg\\
   D-86135 Augsburg, Germany.\\  } 

\end{center}

\begin{abstract}

The quasigeostrophic model is a simplified geophysical fluid model
at asymptotically high rotation rate or at small Rossby number.
We consider the quasigeostrophic equation with dissipation under random
forcing in bounded domains.  We show that global unique solutions
exist for appropriate initial data.  
Unlike the deterministic quasigeostrophic equation whose
well-posedness is well-known, there seems no rigorous result 
on global existence
and uniqueness of the randomly forced quasigeostrophic equation.
Our work provides such a rigorous result on global existence
and uniqueness, under very mild conditions.


\end{abstract} 
 
\newpage
 
{\bf Running head:}

	Random Quasigeostrophic Dynamics 
	
\bigskip
\bigskip

{\bf Author for correspondence:}

	Dr. Jinqiao Duan
	
   Department of Mathematical Sciences  
   
   Clemson University, 
   
   Clemson, South Carolina 29634, USA.
     
   E-mail: duan@math.clemson.edu
   
   Fax: (864) 656-5230
   
   Tel: (864) 656-2730

\newpage
 
\section{Introduction}

The models for geophysical flows are usually too complicated for
analysis.  Simplified partial differential equation models which are
intended to capture the key features of large scale phenomena and
filter out undesired high frequency oscillations in geophysical flows
have been derived at asymptotically high rotation rate or small Rossby
number. An important example of such partial differential equations is
the quasigeostrophic model.

The {\em deterministic} quasigeostrophic equation is
(\cite{Pedlosky},\cite{Benoit})
\begin{equation}
   \Delta \psi_t + J(\psi, \Delta \psi ) + \beta \psi_x
   =\nu \Delta^2 \psi - r \Delta \psi \; ,   
\label{oldeqn}
\end{equation}
where $\psi(x,y,t)$ is the stream function,  
$\beta \geq 0$ is the meridional gradient of the
Coriolis parameter, $\nu >0$ is the viscous dissipation constant,
$r>0$ is the Ekman dissipation constant.
Moreover, $J(f, g)=f_x g_y -f_y g_x$ denotes the Jacobian operator.

The deterministic quasigeostrophic equation (\ref{oldeqn}) has been
derived as an approximation of the rotating shallow water equations by
the conventional asymptotic expansion in small Rossby number
(\cite{Pedlosky}). Schochet (\cite{Schochet}) has recently shown that
the shallow water flows converge to the quasigeostrophic flows in
Sobolev norms in the limit of zero Rossby number (i.e., at
asymptotically high rotation rate), for appropriate initial data. For
related work about the three dimensional baroclinic quasigeostrophic
model, see, for example, \cite{Beale}, \cite{Embid_Majda},
\cite{Holm1} and \cite{Babin1}.

Recently, a few authors have considered the randomly forced
quasigeostrophic equation, in order to incorporate the impact of
uncertain geophysical forces (\cite{Samelson}, 
\cite{Griffa}, \cite{Holloway}, \cite{Muller},
\cite{DelSole-Farrell}).  They studied statistical issues such as
estimating correlation coefficients for the {\em linearized}
quasigeostrophic equation with random forcing. There is also recent
work about the impact of random ocean bottom topography on
quasigeostrophic dynamics (\cite{Klyatskin_Gurarie}).

The randomly forced quasigeostrophic equation takes the form
(\cite{Muller})
\begin{equation}
   \Delta \psi_t + J(\psi, \Delta \psi ) + \beta \psi_x
   =\nu \Delta^2 \psi - r \Delta \psi + \frac{dW}{dt} \; ,
   \label{qg}
\end{equation}
where $W(x,y,t)$ is a space-time Wiener process to be defined
below. There does not seem to exist a mathematically rigorous theory
of quasigeostrophic dynamics under random forcing.  In this paper, we
consider existence and uniqueness of solutions for the nonlinear
quasigeostrophic equation (\ref{qg}) subject to Dirichlet boundary
conditions and appropriate initial data.

\section{Local existence and uniqueness of solution processes}

Introducing $\om (x,y,t) = \Delta \psi(x, y,t)$, the equation (\ref{qg}) can
be written as
\begin{equation}
       \om_t + J(\psi, \om ) + \beta \psi_x
        =\nu \Delta \om - r \om  + \frac{dW}{dt} \; ,
 \label{eqn}
\end{equation}
where $(x, y) \in D$ and $D \subset \R^2$ denotes a   bounded domain
with sufficiently regular boundary.  This equation is supplemented by zero
Dirichlet boundary conditions (\cite{Cessi}) for both $\psi$
and $\om = \Delta \psi$, together with an appropriate initial condition,
i.e., we require
\begin{eqnarray}
 \psi(x,y,t) & = & 0 \quad \mbox{on} \; \partial D \; , \label{BC1} \\
 \om (x,y,t) & = & 0 \quad \mbox{on} \; \partial D \; , \label{BC2} \\
 \om (x,y,0) & = & \om_0(x,y)                    \; . \label{IC}
\end{eqnarray} 
We note that the Poincar\'e inequality holds with these boundary
conditions.

As it stands, (\ref{eqn}) still has to be given a mathematically
precise meaning. This can be done using the framework of stochastic
partial differential equations (\cite{DaPrato}). For this we (formally)
rewrite (\ref{eqn}) in the form
\begin{equation}
    d\om =(\nu \Delta \om - r \om -\beta \psi_x -J(\psi, \om ))dt + dW \; .
    \label{neweqn}
\end{equation}

In the following we use the abbreviations $L^2=L^2(D)$, 
$L^{\infty}=L^{\infty}(D)$,
$H^k_0=H^k_0(D)$, $H^k =H^k(D)$, $ 0 < k < \infty$, for the standard Sobolev
spaces.  Let $<\cdot, \cdot>$ and $\| \cdot \| \equiv \| \cdot \|_2$
denote the standard scalar product and norm in $L^2$, respectively.
Moreover, the norms for $H^k_0$, $L^{\infty}$ are denoted by 
$\| \cdot \|_{H^k}$, $\| \cdot \|_{\infty}$, respectively. 
Due to the Poincar\'e inequality (\cite{Gilbarg-Trudinger}, p.~164), 
$\| \Delta \phi \|$
is an equivalent norm for  $H^2_0$.
 It is well-known that the operator $A = \nu \Delta : L^2
\to L^2$ with domain $D(A) = H^2 \cap H_0^1$ is self-adjoint.  Note
that $A$ generates an analytic semigroup $S(t)$ on $L^2$
(\cite{Pazy}).  The spectrum of $ A$ consists of eigenvalues $0 >
\lambda_1 > \lambda_2 \ge \lambda_3 \ge \ldots$ with corresponding
normalized eigenfunctions $\phi_1$, $\phi_2$, $\ldots$. The set of
these eigenfunctions is complete in $L^2$. For example, for the square
domain $D = (0,1) \times (0,1)$ the eigenvalues are given by
$-\nu(m^2+n^2)\pi^2$ for $m,n \in \N$, and the associated
eigenfunctions are suitable multiples of $\sin(m\pi x) \sin(n\pi y)$.

Now we can define an appropriate class of Wiener processes $W$.  Let
$\beta_k(t)$, $k \in \N$, denote a family of independent real-valued
Brownian motions. Furthermore, choose positive constants $\mu_k$, $k
\in \N$, such that
\begin{displaymath}
  \sum_{k=1}^\infty \frac{\mu_k}{|\lambda_k|^{1-\gamma}} < \infty
\end{displaymath}
for some $0 < \gamma < 1$. Then we consider the Wiener process $W$ defined
by
\begin{eqnarray}
  W(t) := \sum_{k=1}^\infty \sqrt{\mu_k} \beta_k(t) \phi_k \; ,
  \quad t \ge 0 \; .
  \label{Wiener}
\end{eqnarray}
We further assume that  
$$
\kappa (D) = \inf_{0<\rho < diam(D)} \inf_{(x,y) \in D}
		\frac{meas(D \cap B(x,y; \rho))}{\rho^2} > 0,
$$
where $	diam(D) $ is the diameter of $D$ (the least upper bound
of two-point distances in $D$),   $ meas(\cdot)$
denotes the Lebesgue measure, and $B(x,y; \rho)$
is the open disk centered at $(x,y)$ and with radius $\rho$. 
We also assume that
the eigenfunctions $\phi_k$ satisfy
\[
\phi_k \in C_0(\bar{D}), \quad |\phi_k(x,y)| \leq C,
\]
\[
|\p_x \phi_k(x,y)|,  \quad |\p_y \phi_k(x,y)| \leq C \sqrt{|\lambda_k|},
\]
for $(x,y) \in D, k \in \N$, and some constant $C>0$. For the square
domain $D = (0,1) \times (0,1)$, these conditions are
all satisfied.
Then, according to Theorem 5.2.9 in \cite{DaPrato2}, the
stochastic convolution
\begin{eqnarray}
   W_A(t)  =  \int_0^t S(t-s)dW(s) \; , \quad t>0 \; ,
  \label{conv}
\end{eqnarray}
has a continuous version with values in $C_0(D)$,
the Banach space of continuous functions satisfying zero Dirichlet
boundary conditions on $D$.

If we define the nonlinear operator $F$ by $F(\om) = -r \om
-\beta\psi_x -J(\psi, \om )$, then (\ref{neweqn}) can be rewritten as
the abstract evolution equation together with initial condition
\begin{eqnarray}
    d\om & = & (A \om + F(\om) )dt + dW \; ,  
  		\label{QG}     \\
    & & \om (0) = \om_0 
    		\label{newIC}
\end{eqnarray}
or in the mild (integral) form
\begin{eqnarray}
     \om(t) = S(t) \om_0 + \int_0^t S(t-s)F(\om (s)  )ds + W_A(t) \; ,   
  		\label{mild}  
\end{eqnarray}
where the stochastic convolution $W_A(t)$ is defined in (\ref{conv}).

By defining $U(t)=\om(t)-W_A(t)$, we obtain a deterministic mild (integral)
equation
\begin{eqnarray}
     U(t) = S(t) \om_0 + \int_0^t S(t-s)F( U(s)+W_A(s)  )ds \; ,   
  		\label{newmild}  
\end{eqnarray}
or in its differential form
\begin{eqnarray}
    U'(t) & = & A U(t) +  F( U(t) + W_A(t) )  \; ,  
  		\label{newQG}   \\
    & & U(0) = \om_0. \label{newnewIC}
\end{eqnarray}

In the following we prove the local existence of $U(t)$.  We follow
the approach in \cite{DaPrato3} or \cite{DaPrato2}, p.~261.  We first
show that the integral in (\ref{newmild}) makes sense for $U \in
C([0,T];L^2)$. Then, we obtain local existence for (\ref{newmild}) by
the Banach contraction mapping principle in $L^2$.

Note that since $A$ generates an analytic semigroup $S(t)$ on $L^2$ and
has only negative eigenvalues, we have (\cite{Pazy}, p.~74), for $a>0$,
\begin{eqnarray}
    S(t)(-A)^a & = & (-A)^a S(t),  \\
   \| (-A)^a S(t)u \|  & \leq & \frac{c}{t^a} \cdot \|u\|, \\
    \|  S(t)u \|    & \leq &   c  \cdot \|u\|.
   \label{semigroup1}
\end{eqnarray}
Here and hereafter we use $c$ to denote various constants.

We first show that $\int_0^t S(t-s)F( U(s)+W_A(s) )ds$ makes sense for
$U(\cdot)+W_A(\cdot)$ (and thus $U(\cdot)$) in
$C([0,T];L^2)$. Recalling that $\om = U + W_A$, this follows from the
following lemma.

\begin{lemma} \label{f_global}
  Define the mapping $\cF : C( [0,T]; H_0^1) \to C( [0,T]; L^2)$ by
  \begin{displaymath}
    (\cF(\omega))(t) := \int_0^t S(t-s) F(\omega(s)) ds \; ,
    \quad t \in [0,T] \; ,
    \quad \omega \in C( [0,T]; L^2) \; .
  \end{displaymath}
  Then $\cF$ is continuous, and it can be extended to a continuous
  mapping from the space $C( [0,T]; L^2)$ to $C([0,T]; L^2)$.
  Furthermore, the image of the extended mapping $\cF$ is contained
  in $C([0,T],H^a(D))$ for $0 \le a < \frac12$.
\end{lemma}

\noindent
{\bf Proof:} The continuity of $\cF : C( [0,T]; H_0^1) \to C( [0,T]; L^2)$
is obvious. As for extending the domain of $\cF$ let $\omega, \bar\omega
\in C( [0,T]; L^2)$ be arbitrary. Using the abbreviations $\om = \Delta \psi$
and $\bar\om =\Delta \bar\psi$ we get
\begin{eqnarray*}
F(\om)- F(\bar\om) 
  & = &  r (\bar\om-\om)+ \beta (\bar\psi-\psi)_x +
	\psi_x (\bar\om-\om)_y + (\bar\psi -\psi)_x \bar\om_y   \\
  & & \quad - \psi_y (\bar\om-\om)_x +(\psi-\bar\psi)_y \bar\om_x . 
\end{eqnarray*}
Let $a \in [0, 1)$, and consider an arbitrary $\phi \in D( (-A)^a )$.
Then the above identity implies
\begin{eqnarray*}
  I & = & \left\langle (-A)^a \phi, \int_0^t S(t-s) F(\om(s))ds -
    \int_0^t S(t-s) F(\bar\om(s))ds \right\rangle \\
  & = & \int_0^t \langle S(t-s)(-A)^a \phi,
    [ r (\bar\om-\om)+ \beta (\bar\psi-\psi)_x +
    \psi_x (\bar\om-\om)_y \\
  & & \qquad\quad + (\bar\psi - \psi)_x \bar\om_y 
    - \psi_y (\bar\om-\om)_x +(\psi-\bar\psi)_y \bar\om_x](s)
    \rangle ds \\
& \equiv  &  \int_0^t (I_1 + I_2 + I_3 + I_4 + I_5 + I_6) ds,
\end{eqnarray*}
where
\begin{eqnarray*} 
I_1 & = & <S(t-s)(-A)^a \phi , [r (\bar\om-\om)](s) >, \\
I_2 & = & < S(t-s)(-A)^a \phi , [\beta (\bar\psi-\psi)_x](s) >, \\
I_3 & = & < S(t-s)(-A)^a \phi , [\psi_x (\bar\om-\om)_y](s) >, \\
I_4 & = & <S(t-s)(-A)^a \phi, [(\bar\psi - \psi)_x \bar\om_y](s)>, \\
I_5 & = & <S(t-s)(-A)^a \phi, [- \psi_y (\bar\om-\om)_x](s)>, \\
I_6 & = & <S(t-s)(-A)^a \phi, [(\psi-\bar\psi)_y \bar\om_x](s)>.
\end{eqnarray*}
Now we estimate $|I_k|$, $k=1,\ldots,6$, one by one, thereby omitting
the argument $s$.
\begin{eqnarray*}
|I_1| & = & |<S(t-s)(-A)^a \phi , r (\bar\om-\om)>|       \\
& \leq & r \|S(t-s)(-A)^a \phi \| \cdot \|\bar\om-\om\|     \\
& \leq & rc (t-s)^{-a} \cdot \|\phi\| \cdot \|\bar\om-\om\|,
\label{estimate1}
\end{eqnarray*} 
\begin{eqnarray*}
|I_2| & = & |< S(t-s)(-A)^a \phi , \beta (\bar\psi-\psi)_x >| \\
& \leq &  \beta \| S(t-s)(-A)^a \phi \| \cdot \|(\bar\psi-\psi)_x\|    \\
& \leq & \beta c (t-s)^{-a} \cdot \|\phi\| \cdot \|\bar\om-\om\|,
\label{estimate2}
\end{eqnarray*}
where we have used the Poincar\'e inequality
(\cite{Gilbarg-Trudinger}, p.~164) on $(\bar\psi-\psi)_x$ which has
zero mean. Using the Cauchy-Schwarz inequality $|<u,v>| \leq \|u\|
\|v\|$ we also obtain
\begin{eqnarray*}
|I_3| & = & |< S(t-s)(-A)^a \phi , \psi_x(\bar\om-\om)_y>|  \\
& = & |<D_y [(S(t-s)(-A)^a \phi) \psi_x], \bar\om - \om>|  \\
& \leq & |<D_y (S(t-s)(-A)^a \phi )\psi_x, \bar\om - \om> | \\
& & \qquad + |< (S(t-s)(-A)^a \phi) \psi_{xy}, \bar\om - \om>|  \\
& \leq & \|D_y (S(t-s)(-A)^a \phi)\psi_x\| \cdot \|\bar\om - \om\| \\
& & \qquad +
    \|S(t-s)(-A)^a \phi\|_{\infty} \cdot \|\psi_{xy}\| \cdot
    \|\bar\om - \om\|.
\end{eqnarray*}
As for estimating $\|D_y (S(t-s)(-A)^a \phi)\psi_x\|$ we get
\begin{eqnarray*}
\|D_y (S(t-s)(-A)^a \phi)\psi_x\| & \le &
    \|D_y (S(t-s)(-A)^a \phi)\|_{\frac14} \cdot \|\psi_x\|_{\frac14} \\
& \leq & c \|D_y (S(t-s)(-A)^a \phi)\|_{H^{\frac12}} \cdot
    \|\psi_x\|_{H^1} \\
& \leq & c \|S(t-s)(-A)^a \phi\|_{H^{\frac32}} \cdot \|\om\| \\
& \leq & c \|A^{\frac34+\rho} S(t-s)(-A)^a \phi\| \cdot \|\om\| \\
& \leq & c(t-s)^{-(\frac34+\rho+a)} \cdot \|\phi\| \cdot \|\om\| ,
\end{eqnarray*}
where we have used the inequality $\|uv\| \leq \|u\|_{\frac14}
\|v\|_{\frac14}$, the continuity of the mapping $D_y : H^{\frac32} \to
H^{\frac12}$ (\cite{Temam}, p.~56), the embedding $D( (-A)^{\frac34 +
\rho} ) \hookrightarrow H^{\frac32}$ for arbitrary $\rho > 0$
(\cite{Pazy}, p.~243), and the facts that $H^{\frac12}$ and $H^1$ are
embedded in $L^4$ (\cite{Adams}, p.~217). Furthermore,
\begin{eqnarray*}
\|S(t-s)(-A)^a \phi\|_{\infty} & \leq &
    c \|S(t-s)(-A)^a \phi\|_{H^{1 + \rho}} \\
& \leq & c \|(-A)^{\frac12 + \rho} S(t-s) (-A)^a \phi\| \\
& \leq & c(t-s)^{-(\frac12+\rho+a)} \cdot \|\phi\|,
\end{eqnarray*}
due to the smoothing property of the semigroup $S$ and the embeddings
$D( (-A)^{\frac12 + \rho} ) \hookrightarrow H^{1 + \rho}
\hookrightarrow L^{\infty}$ for arbitrary $\rho > 0$ (\cite{Pazy},
pp.~208, 243). Altogether we get
\begin{displaymath}
|I_3| \leq [c(t-s)^{-(\frac34+\rho+a)} + c(t-s)^{-(\frac12+\rho+a)}]
    \cdot \|\phi\| \cdot \|\om\| \cdot \|\om - \bar\om\|.
\end{displaymath}
Similarly, the following estimates can be obtained:
\begin{eqnarray*}
  |I_4| & \le & [c(t-s)^{-(\frac34+\rho+a)} +
    c(t-s)^{-(\frac12+\rho+a)}] \cdot \|\phi\| \cdot
    \|\bar\om\| \cdot \|\om - \bar\om\| \; , \\
  |I_5| & \le & [c(t-s)^{-(\frac34+\rho+a)} +
    c(t-s)^{-(\frac12+\rho+a)}] \cdot \|\phi\| \cdot
    \|\om\| \cdot \|\om - \bar\om\| \; , \\
  |I_6| & \le & [c(t-s)^{-(\frac34+\rho+a)} +
    c(t-s)^{-(\frac12+\rho+a)}] \cdot \|\phi\| \cdot
    \|\bar\om\| \cdot \|\om - \bar\om\| \; .
\end{eqnarray*}
Thus, we have
\begin{eqnarray*}
|I|  
& \leq &  \int_0^t (|I_1| + |I_2| + |I_3| + |I_4| + |I_5| + |I_6|) ds   \\
& \leq &  \frac{rc + \beta c}{1-a} \cdot t^{1-a} \cdot
    \|\phi\| \cdot \sup\limits_{0\leq s \leq t} \|\om(s)-\bar\om(s)\|   \\
& & \quad + \left( \frac{8c}{1-4\rho-4a} \cdot t^{\frac14-\rho-a} +
    \frac{4c}{1-2\rho-2a} \cdot t^{\frac12-\rho-a} \right) \\
& & \qquad \cdot \, \|\phi\| \cdot \sup\limits_{0\leq s \leq t}
  ( \|\om(s)\| + \|\bar\om(s)\| ) \cdot
  \sup\limits_{0\leq s \leq t} \|\om(s)-\bar\om(s)\|,
\end{eqnarray*}
provided the positive constants $a$ and $\rho$ satisfy $0 < \rho + a <
\frac14$. This finally implies
\begin{displaymath}
  \int_0^t S(t-s) F(\om(s)) ds - \int_0^t S(t-s) F(\bar\om(s)) ds
  \in D( (-A)^a)
  \quad\mbox{for}\quad
  0 \leq a < \frac14,
\end{displaymath}
and
\begin{eqnarray*} 
& & \left\| (-A)^a \left(\int_0^t S(t-s) F(\om (s)) ds -
  \int_0^t S(t-s) F(\bar\om (s)) ds \right) \right\| \\
& \leq & \frac{rc + \beta c}{1-a} \cdot t^{1-a} \cdot
    \sup\limits_{0\leq s \leq t} \|\om(s)-\bar\om(s)\|   \\
& & \quad + \left( \frac{8c}{1-4\rho-4a} \cdot t^{\frac14-\rho-a} +
    \frac{4c}{1-2\rho-2a} \cdot t^{\frac12-\rho-a} \right) \\
& & \qquad \cdot \sup\limits_{0\leq s \leq t}
  ( \|\om(s)\| + \|\bar\om(s)\| ) \cdot
  \sup\limits_{0\leq s \leq t} \|\om(s)-\bar\om(s)\|.
\end{eqnarray*}
Especially for $a=0$ we obtain
\begin{eqnarray*} 
& & \left\| \int_0^t S(t-s) F(\om (s)) ds -
  \int_0^t S(t-s) F(\bar\om (s)) ds \right\| \\
& \leq & (rc + \beta c) \cdot t \cdot
    \sup\limits_{0\leq s \leq t} \|\om(s)-\bar\om(s)\|   \\
& & \quad + \left( \frac{8c}{1-4\rho} \cdot t^{\frac14-\rho} +
    \frac{4c}{1-2\rho} \cdot t^{\frac12-\rho} \right) \\
& & \qquad \cdot \sup\limits_{0\leq s \leq t}
  ( \|\om(s)\| + \|\bar\om(s)\| ) \cdot
  \sup\limits_{0\leq s \leq t} \|\om(s)-\bar\om(s)\|,
\end{eqnarray*}
for every $0 < \rho < \frac14$. This completes the proof of the
lemma.
\hspace*{\fill}$\Box$

\bigskip
We conclude from the above lemma that $\int_0^t S(t-s) F(U(s) +
W_A(s)) ds$, considered as a mapping with argument $U(\cdot)$, can be
extended to a bounded map from $C([0, T], L^2(D))$ into itself.

Now we can follow \cite{Pazy}, p.~196 or \cite{DaPrato}, p.~201, to
obtain that (\ref{newmild}) has a unique local solution $U(t)$, or
(\ref{QG}), (\ref{newIC}) has a unique local solution $\om(x, y, t)$,
on $[0, \tau)$, by the Banach contraction mapping principle. The
solution $\om(x,y,t)$ is in   $C([0, \tau]; L^2(D))$, as well as
in $C((0, \tau]; H^a(D))$, for arbitrary $0 \le a < \frac12$.

\section{Global solution processes}

In this section, we show that the  solution $U(t)$ is a priori bounded,
in $L^2(D)$-norm, on any finite interval $[0, T]$. This implies that the
local solution $U(t)$, and thus $\om(x, y, t)$ is actually global in time.

We consider (\ref{newQG}), (\ref{newnewIC}) with $W_A$ replaced by a
regular function $V$ from the space $C([0, T]; H^3_0(D))$ and $\om_0 $
in $D(A)$
\begin{eqnarray}
    U'(t) & = & A U(t) +  F( U(t) + V(t) )  \; ,   \label{regulareqn}  \\
    & & U (0) = \om_0, \nonumber
\end{eqnarray}
where, we denote, $U=\Delta u, V=\Delta v$.  More specifically,
(\ref{regulareqn}) is
\begin{eqnarray}
    U' = \nu \Delta U-r(U+V) -\beta (u+v)_x -J(u+v, U+V).
    \label{regulareqn1} 
\end{eqnarray}
Due to the smoothing effect of the sectorial operator $A$, and the
fact that $F$ is locally Lipschitz in $U$ from $H_0^{m+1} \cap
H^{2m+2}$ to $H_0^{m} \cap H^{2m}$ for $m=0, 1, 2$, we
conclude that the solution $U$ of (\ref{regulareqn}) is in $H_0^{k}
\cap H^{2k}$ for   $k=0, 1, 2$, and hence $U$ is a strong
solution (if $V$ is smoother, the solution $U$ is also
smoother); see \cite{Henry}, p.~73.
  
We now estimate the norm $\| U(t) \|$.

Multiplying (\ref{regulareqn1}) by $U$ and integrating over $D$, we get
\begin{eqnarray}
\frac12 \frac{d}{dt}\|U\|^2 
& = & -\nu \int_D |\nabla U|^2 -r\int_D (U + V)U       \nonumber  \\
& & \quad - \beta \int_D (u_x +v_x)U  - \int_D J(u+v, U+V)U   \nonumber  \\
& = &  -\nu \int_D |\nabla U|^2 -r\int_D (U^2+UV)       \nonumber  \\ 
& & \quad - \beta \int_D (u_xU+v_xU) - \int_D J(u,  V)U - \int_D J(v, V)U 
					\nonumber  \\
& = &  -\nu  \|\nabla U \|^2 -r\int_D (U^2+UV)
    - \beta \int_D (u_xU+v_xU)		\nonumber  \\
& & \quad + \int_D (-u_xV_yU +u_yV_xU -v_xV_yU + v_yV_xU),		
 	\label{term1} 
\end{eqnarray} 
where we have used the fact that $\int_D J(u, U)U = \int_D J(v, U)U =0$
via integration by parts; see also  \cite{Holm1}.
We estimate the right hand side of (\ref{term1}) term by term.
\begin{eqnarray}
  -r\int_D (U + V)U  
& \leq &  r (1+c\|V\|_{\infty}) \|U\|^2  ,
\label{term2}
\end{eqnarray}
\begin{eqnarray}       
 - \beta \int_D (u_xU+v_xU)  
& \leq & \beta \int_D \frac12[u_x^2 +U^2 +v_x^2 + U^2]  \nonumber  \\
& \leq & \beta c (\|U\|^2 + \|V\|^2)  \nonumber  \\
& \leq & \beta c (\|U\|^2 + \|V\|_{\infty}^2),
\label{term3}
\end{eqnarray}
where we have used the Poincar\'e inequality on $u_x$, $v_x$, which
have zero mean.
\begin{eqnarray}  
 \int_D (-u_xV_yU) 
 & = &  \int_D (u_xU)_y V    		    \nonumber  \\	
 & = & \int_D (u_{xy}UV +u_xU_y V)  	    \nonumber  \\
 & \leq & \|V\|_{\infty} \int_D |u_{xy}U| +\int_D (|u_x|\cdot\|V\|_{\infty})
    |U_y|  \nonumber  \\
 & \leq & \|V\|_{\infty} \int_D \frac12 (u_{xy}^2 +U^2)   
 	+\int_D \left(\frac1{2\e} u_x^2 \|V\|_{\infty}^2 +\frac{\e}2 U_y^2
        \right)  
 						 \nonumber  \\
 & \leq &  \frac{c}{2} \|V\|_{\infty} \left(1+\frac1{\e}\|V\|_{\infty}\right)
      \|U\|^2
 	 + \frac{\e}2 \|\nabla U\|^2  ,
 \label{term4}
\end{eqnarray} 
since $\int_D u_{xy}^2$ is bounded by $c \int_D (\Delta u)^2 = c
\int_D U^2$. We have also used the Young inequality (\cite{Temam}) to
get that $(|u_x|\cdot\|V\|_{\infty}) |U_y| \leq \frac1{2\e}
u_x^2\|V\|_{\infty}^2 + \frac{\e}2 U_y^2$, for any positive real
number $\e >0$.  Similarly, we have
\begin{eqnarray}  
\int_D  u_yV_xU  
& \leq &   \frac{c}{2} \|V\|_{\infty} \left(1+\frac1{\e}\|V\|_{\infty}\right)
    \|U\|^2
 	 + \frac{\e}2 \|\nabla U\|^2 ,    
  			    \label{term5}   \nonumber  \\
\int_D (-v_xV_yU )
& = &  \int_D (v_xU)_y V 	      \nonumber  \\
& = & \int_D (v_{xy}UV + v_xU_y V)    \nonumber  \\
& \leq & \|V\|_{\infty} \int_D \frac12 (v_{xy}^2 +U^2)   
 	+\int_D \left(\frac1{2\e} v_x^2 \|V\|_{\infty}^2 +
        \frac{\e}2 U_y^2 \right)
  					\nonumber  \\
& \leq &  \frac12 \|V\|_{\infty} \|U\|^2 +c\|V\|_{\infty} \|V\|^2
	+\frac{c}{\e} \|V\|_{\infty}^2\|V\|^2+ \frac{\e}2 \|\nabla U\|^2 
 			\nonumber  \\
& \leq &  \frac12 \|V\|_{\infty} \|U\|^2 + c\|V\|_{\infty}^3
		+ \frac{c}{\e} \|V\|_{\infty}^4 + \frac{\e}2 \|\nabla U\|^2,
    				\label{term6}   \nonumber  \\
\int_D v_yV_xU  
 & \leq &  \frac12 \|V\|_{\infty} \|U\|^2+ c\|V\|_{\infty}^3
		+ \frac{c}{\e} \|V\|_{\infty}^4 + \frac{\e}2 \|\nabla U\|^2.
    				\label{term7}   				
\end{eqnarray}
Putting (\ref{term2}) - (\ref{term7}) into  (\ref{term1}), we get
\begin{eqnarray}
\frac12 \frac{d}{dt}\|U\|^2
&\leq &   (-\nu +  2 \e  )\cdot \|\nabla U\|^2   \nonumber  \\
&  & +[ r (1+c\|V\|_{\infty}) + \b c 
 + c \|V\|_{\infty} (1+\frac1{\e}\|V\|_{\infty}) +
    \|V\|_{\infty}]
    \cdot \|U\|^2	\nonumber  \\
&  &  +  \b c \|V\|^2_{\infty} + 2c\|V\|_{\infty}^3 
	+ \frac2{\e} c \|V\|_{\infty}^4 .
\end{eqnarray}
Taking $\e = \frac{\nu}2$, we finally obtain
\begin{eqnarray}
 \frac{d}{dt}\|U(t)\|^2
& \leq &  A(t) \cdot \|U(t)\|^2 + B(t),
\end{eqnarray}
where
\begin{eqnarray}
A(t) & = & 2 \left[ r (1+c\|V\|_{\infty}) + \b c 
  + c \|V\|_{\infty} \left(1+ \frac2{\nu}\|V\|_{\infty}\right)
  + \|V\|_{\infty} \right] > 0,  \nonumber  \\
B(t) & = & 2\b c \|V\|^2_{\infty} +  4 c \|V\|_{\infty}^3
	+  \frac8{\nu}c \|V\|_{\infty}^4 > 0.
\end{eqnarray}
Hence by the Gronwall inequality (\cite{Temam}), we obtain
\begin{equation}
  \|U(t)\|^2  \leq \|\om_0\|^2 e^{\int_0^t A(s) ds}
  +\int_0^t B(s) e^{\int_s^t A(\tau)d\tau} ds, \; 0<t<T.
  \label{bounded}
\end{equation}
Note that $H^3_0(D)$ is embedded in $C_0(D)$, the trajectories of 
$W_A(t)$ can be uniformly approximated, on any finite interval $[0, T]$,
by functions $V$ in $C([0, T]; H^3_0(D))$, and
$D(A)$ is dense in $L^2(D)$. Thus the boundedness estimate
(\ref{bounded}) is true for any local solution $U(t)$ of (\ref{newQG}).
This shows that the unique (local) solution does not blow 
up on any finite intervals. 

We thus have the following theorem.

\begin{theorem}
  For every initial condition $\om_0(x,y) \in L^2(D)$, there exists a
  unique global mild solution $\om(x,y,t)$ of the quasigeostrophic model
  (\ref{eqn}), (\ref{BC1}), (\ref{BC2}), and (\ref{IC}). This solution
  is contained in the space $C([0, T]; L^2(D))$ for every $T>0$, as well
  as in $C((0, T]; H^a(D))$ for all $0 \le a < \frac12$ and $T>0$.
\end{theorem}

\section{Discussions}

There has been recent work on geophysical problems modeled
by the randomly forced quasigeostrophic equation
(e.g., \cite{Samelson}, \cite{Griffa}, \cite{Holloway},
\cite{Muller}, \cite{DelSole-Farrell}).
Unlike the deterministic quasigeostrophic equation whose
well-posedness is well-known (e.g., \cite{Bennett}, \cite{Schochet},
\cite{Beale}), there seems no   rigorous results on global existence
and uniqueness of the randomly forced quasigeostrophic equation.
Our work provides such a rigorous result on global existence
and uniqueness, under very mild conditions.

\end{document}